\documentclass[11pt]{amsart}
\usepackage{url}
\newcommand{\cf}{{\operatorname{cf}}}
\newcommand{\pcf}{{\operatorname{Pcf}}}
\newcommand{\REG}{\operatorname{Reg}}
\newcommand{\cov}{\operatorname{Cov}}
\newtheorem{theorem}{Theorem}
\newtheorem{definition}[theorem]{Definition}

\begin{document}
\noindent
{\small Topology Atlas Invited Contributions {\bf 6} no.~4 (2001) 4 pp.}
\vspace{\baselineskip}
\title{PCF Theory}
\author{Menachem Kojman}
\address{Department of Mathematics\\
Ben Gurion University of the Negev\\
Beer Sheva, Israel}
\email{kojman@math.bgu.ac.il}
\maketitle 
\section{PCF theory and singular cardinals}

The abbreviation ``PCF'' stands for ``Possible Cofinalities''. 
PCF theory was invented by Saharon Shelah to prove upper bounds on 
exponents of singular cardinals.

The starting point of PCF theory is in the realization that the usual 
exponent function is too coarse for measuring the power set of singular 
cardinals.

Consider the cardinal $\aleph_\omega$, which is the smallest singular
cardinal and has countable cofinality. 
The usual exponent $(\aleph_\omega)^{\aleph_0}$ measures the total number
of countable subsets of $\aleph_\omega$ and is larger than $\aleph_\omega$
itself by straightforward diagonalization. 
The exponent $(\aleph_\omega)^{\aleph_0}$ clearly satisfies 
$(\aleph_\omega)^{\aleph_0} \ge 2^{\aleph_0}$ and therefore has no upper
bound in the list $\{\aleph_\alpha : \alpha\in \mathrm{On}\}$ of cardinal 
numbers (because $2^{\aleph_0}$ itself has none).

The main conceptual change that PCF theory has generated lies in the fact 
that the number of countable subsets of $\aleph_\omega$ which is required 
to \emph{cover} all countable subsets of $\aleph_\omega$ is always bounded 
by $\aleph_{\omega_4}$.

\begin{definition} $\cov(\aleph_\omega,\omega)$ is the least number of 
countable subsets of $\aleph_\omega$ required to cover all countable 
subsets of $\aleph_\omega$. 
Equivalently, $\cov(\aleph_\omega,\omega)$ is the \emph{cofinality} of the
partially ordered set $([\aleph_\omega]^{\aleph_0},\subseteq)$.
\end{definition}

Fixing a covering collection $\mathcal{F}\subseteq
[\aleph_\omega]^{\aleph_0}$ of cardinality $\cov(\aleph_\omega,\omega)$
and replacing each $A\in \mathcal{F}$ by $[A]^{\aleph_0}$ one gets the
following equality:
$$
(\aleph_\omega)^{\aleph_0}=\cov(\aleph_\omega,\omega)\times 2^{\aleph_0}.
$$

\begin{theorem}[Shelah] $\cov(\aleph_\omega,\omega)<\aleph_{\omega_4}$.
\end{theorem}

The factor $2^{\aleph_0}$ in $(\aleph_{\omega})^{\aleph_0}$ is not 
bounded; the factor $\cov(\aleph_\omega,\omega)$ is. 
The number $\cov(\aleph_\omega,\omega)$ cannot be changed by ccc forcing
(since every new countable set will be contained in an old countable
set) and, furthermore, to make $\cov(\aleph_\omega,\omega)$ greater than 
$\aleph_{\omega+1}$, its minimal possible value, requires the consistency 
of large cardinals.

From this point of view, the power $2^{\aleph_0}$ of the regular cardinal
$\aleph_0$ generates some ``noise'' in the collection 
$[\aleph_\omega]^{\aleph_0}$, whose important properties are captured by
the robust factor $\cov(\aleph_\omega,\omega)$. 

\section{Reduced products and scales}

The way PCF theory computes covering numbers of singular cardinals is by
reducing them to the algebra of reduced products (modulo ideals) of sets
of regular cardinals $A\subseteq \REG$ with the condition $\min A > |A|$. 
For such a set $A$, the set of \emph{possible cofinalities} is defined by:
$$
\pcf A 
= \{ \cf\prod A/U: U\subseteq \mathcal{P}(A) \text { an ultrafilter} \}.
$$

Reduced products of the form $\prod A/I$ where $A\subseteq \REG$ with 
$\min A > A$ modulo an ideal $I$ over $A$ behave very differently from 
products of $\lambda$ copies of a regular cardinal $\lambda$ modulo some 
ideal (e.g.\ $(\omega^\omega,<^*)$). 
The \emph{PCF theorem} is the fundamental structure theorem for products
$\prod A/I$ for such $A\subseteq \REG$. 
It is of course worth pointing out the fact that there \emph{exists} a 
structure theorem at all in this case!

The PCF theorem states the existence of \emph{PCF generators}. 
It shows how each of the possible cofinalities in $\pcf A$ can be 
represented as a scale:

\begin{theorem}[Shelah] 
For every $A\subseteq \REG$ with $\min A > A$ and every $\lambda \in \pcf A$ 
there exists a set $B_\lambda \subseteq A$ called the 
\emph{generator}
of $\lambda$ so that: the product $\prod B_\lambda$ modulo the ideal 
$\langle B_\theta : {\theta \in \pcf A} \wedge {\theta < \lambda}\rangle$ 
has an increasing and cofinal sequence of length $\lambda$.
\end{theorem}

In a less technical formulation: there is a sequence of ideals over $A$,
each generated over the union of previous ones by a single set 
$B_\lambda$, $\lambda \in \pcf A$, so that the cardinal $\lambda \in \pcf 
A$ is represented as the \emph{true cofinality} of the product of 
$B_\lambda$ over the ideal $J_{<\lambda}$ generated by all smaller 
generators. 
A product has true cofinality if it contains an increasing and cofinal
sequence (which is, e.g.\ the situation in $(\omega^\omega,<^*)$ when 
${\mathfrak b} = {\mathfrak d}$).

The witness to the true cofinality of $\prod B_\lambda/J_{<\lambda}$ is a 
sequence
$\overline{f}^\lambda 
= \langle f_\alpha : \alpha < \lambda\rangle$ 
which is $<_{J_{<\lambda}}$-increasing and cofinal in 
$\prod B_{\lambda}/J_{<\lambda}$. 
Such a sequence is called a \emph{ $\lambda$-scale}. 
Among the first important discoveries of PCF theory was the fact that
there is always an $\aleph_{\omega+1}$-scale.

\begin{theorem}[Shelah]\label{scale} 
There exists an infinite set $B \subseteq \omega$ and a sequence 
$\overline{f} = \langle f_\alpha : \alpha < \aleph_{\omega+1} \rangle 
\subseteq \prod_{n\in B}\omega_n$ 
so that:
\begin{enumerate}
\item $\alpha < \beta < \aleph_{\omega+1} \Rightarrow f_\alpha <^* f_\beta$;
\item for all $f\in \prod_{n\in B}\omega_n$ there exists
$\alpha < \aleph_{\omega+1}$ for which $f<^*f_\alpha$.
\end{enumerate}
\end{theorem}

The relation $f<^* g$ means: $f(n)$ is strictly smaller than $g(n)$ in all
but finitely many $n$s.

\begin{theorem}[Shelah] 
$\pcf \{\aleph_n : n\}$ is an interval of regular cardinals with a last 
element $\max\pcf\{\aleph_n : n < \omega\}$ and 
$\max\pcf\{\aleph_n : n < \omega\}$ is equal to 
$\cov(\aleph_\omega,\omega)$.
\end{theorem}

\section{PCF theory and topology}

As every other important development in set theory, also PCF theory has
consequences in topology. The theory is useful in constructing new spaces
and understanding old space whose properties are related to singular
cardinals. 
More importantly, the new way of seeing things which is suggested by PCF
theory is helpful in phrasing new theorems.

A good example to an application of PCF theory in topology is related to 
M.~E.~Rudin's Dowker space $X^R$, whose cardinality is 
$(\aleph_\omega)^{\aleph_0}$. 
Recall the definition of $X^R$:
$$
X^R
= \{ f \in \prod_{n>1}\omega_n+1 : 
(\exists m)(\forall n > 1)(\aleph_0 < \cf f(n) < \aleph_m) \}.
$$

The topology on $X^R$ is the box topology.

It turns out that neither of the two properties of $X^R$---collectionwise 
normality and the absence of countable paracompactness---has anything to 
do with the factor $2^{\aleph_0}$ in $(\aleph_\omega)^{\aleph_0}$. 
The fact that $X^R$ is Dowker is related only to the structure of PCF 
scales in $\prod _{n>1}\omega_n$.

Fix a continuous $\aleph_{\omega+1}$-scale 
$\overline{f} = 
\langle f_\alpha : \alpha < \aleph_{\omega+1} \rangle$ 
by Theorem~\ref{scale}. 
``Continuous'' means that whenever 
$\langle f_\alpha : \alpha < \delta \rangle$ 
for $\delta < \aleph_{\omega+1}$ with $\cf\delta > \aleph_0$ has a least 
upper bound, then $f_\delta$ itself is such a least upper bound. 
A continuous scale is gotten easily from any scale.

Define the following subspace of $X^R$:
$$
X = 
\{f\in X^R : (\exists \alpha < \aleph_{\omega+1})(f =^* f_\alpha)\}.
$$

The space $X$ has cardinality $\aleph_{\omega+1}$, local character 
$\aleph_\omega$ and weight $\aleph_{\omega+1}$ even when the respective 
characteristics of $X^R$ are much larger than $\aleph_{\omega+1}$. 
But just like $X^R$, $X$ is collectionwise normal and not countably 
paracompact. Thus:

\begin{theorem}[Kojman, Shelah \cite{Dowker}] 
There exists a Dowker space of cardinality and weight $\aleph_{\omega+1}$.
\end{theorem}

Another example of an application of PCF techniques to topology is the 
computation of the Baire number of the space of all uniform ultrafilters
over a singular cardinal of countable cofinality.

\begin{theorem}[Kojman, Shelah \cite{fallen}] 
If $\mu$ is the a singular cardinal of countable cofinality, then exactly 
$\aleph_2$ nowhere-dense sets are required to cover the space of all 
uniform ultrafilters over $\mu$. An ultrafilter over $\mu$ is 
\emph{uniform} if it does not contain a set of cardinality $<\mu$.
\end{theorem}

\end{document}